\documentclass[12pt]{article}

\usepackage{amsfonts,graphicx,euscript}

\newtheorem{lemma}{Lemma}
\newtheorem{theorem}{Theorem}

\title{Semidefinite geometry of the numerical range}

\author{Didier Henrion$^{1,2}$}
\date{\today}

\begin{document}

\maketitle

\footnotetext[1]{LAAS-CNRS, University of Toulouse, France}
\footnotetext[2]{Faculty of Electrical Engineering,
Czech Technical University in Prague, Czech Republic}
\addtocounter{footnote}{2}

\begin{abstract}
The numerical range of a matrix is studied geometrically
via the cone of positive semidefinite matrices (or semidefinite cone
for short). In particular
it is shown that the feasible set of a two-dimensional linear
matrix inequality (LMI), an affine section of the semidefinite cone,
is always dual to the numerical range of a matrix, which is
therefore an affine projection of the semidefinite cone. Both primal and dual
sets can also be viewed as convex hulls of explicit algebraic plane
curve components. Several numerical examples illustrate this interplay
between algebra, geometry and semidefinite programming duality.
Finally, these techniques are used to revisit a theorem in statistics
on the independence of quadratic forms in a normally distributed vector.
\end{abstract}

\begin{center}
\end{center}

\section{Notations and definitions}

The numerical range of a matrix $A \in {\mathbb C}^{n\times n}$
is defined as
\begin{equation}\label{wa}
{\mathcal W}(A) = \{ w^* A w \in {\mathbb C} \: :\: w \in {\mathbb C}^n, \: w^*w = 1\}.
\end{equation}
It is a convex closed set of the complex plane which
contains the spectrum of $A$. It is also called the field of values,
see \cite[Chapter 1]{gustafson} and \cite[Chapter 1]{horn}
for elementary introductions. Matlab functions for visualizing numerical
ranges are freely available from \cite{cowen} and \cite{higham}.

Let
\begin{equation}\label{ai}
A_0 = I_n, \quad A_1 = \frac{A+A^*}{2}, \quad A_2 = \frac{A-A^*}{2i}
\end{equation}
with $I_n$ denoting the identity matrix of size $n$ and
$i$ denoting the imaginary unit.
Define
\begin{equation}\label{fa}
{\mathcal F}(A) = \{ y \in {\mathbb P}^2_+ \: :\: F(y) = y_0 A_0 + y_1 A_1 + y_2 A_2 \succeq 0\}
\end{equation}
with $\succeq 0$ meaning positive semidefinite (since the $A_i$ are Hermitian matrices, $F(y)$ has
real eigenvalues for all $y$) and ${\mathbb P}^2_+$ denoting
the oriented projective real plane (a model of the projective plane
where the signs of homogeneous coordinates are significant, and which
allows orientation, ordering and separation tests such as inequalities,
see \cite{stolfi}
for more details).
Set ${\mathcal F}(A)$ is a
linear section of the cone of positive semidefinite matrices (or semidefinite
cone for short), see \cite[Chapter 4]{bental}.
Inequality $F(y) \succeq 0$ is called a linear
matrix inequality (LMI). In the complex plane $\mathbb C$, or equivalently, in the
affine real plane ${\mathbb R}^2$,
set ${\mathcal F}(A)$ is a convex set including the origin, an affine section
of the semidefinite cone.

Let
\[
p(y) = \det(y_0 A_0 + y_1 A_1 + y_2 A_2)
\]
be a trivariate form of degree $n$ defining the algebraic plane curve
\begin{equation}\label{curvep}
{\mathcal P} = \{ y \in {\mathbb P}^2_+ \: :\: p(y) = 0\}.
\end{equation}
Let
\begin{equation}\label{curveq}
{\mathcal Q} = \{ x \in {\mathbb P}^2_+ \: :\: q(x) = 0\}
\end{equation}
be the algebraic plane curve dual to $\mathcal P$,
in the sense that we associate to each point $y \in \mathcal P$
a point $x \in \mathcal Q$ of projective coordinates
$x = (\partial p(y)/\partial y_0,\:\partial p(y)/\partial y_1,\:
\partial p(y)/\partial y_2)$. Geometrically, a point in $\mathcal Q$
corresponds to a tangent at the corresponding point in $\mathcal Q$,
and conversely, see \cite[Section V.8]{walker} and \cite[Section 1.1]{gelfand}
for elementary properties of dual curves.

Let $\mathbb V$ denote a vector space equipped
with inner product $\langle .,.\rangle$. If $x$ and $y$ are vectors
then $\langle x, y \rangle = x^* y$. If $X$ and $Y$ are symmetric
matrices, then $\langle X, Y \rangle = \mathrm{trace}(X^*Y)$.
Given a set $\mathcal K$ in $\mathbb V$,
its dual set consists of all
linear maps from $\mathcal K$ to non-negative elements in $\mathbb R$, namely
\[
{\mathcal K}^* = \{y \in {\mathbb V} \: :\: \langle x,y \rangle \geq 0,
\: x \in {\mathcal K}\}.
\]
Finally, the convex hull of a set $\mathcal K$, denoted $\mathrm{conv}\:{\mathcal K}$,
is the set of all convex combinations of elements in $\mathcal K$.
 
\section{Semidefinite duality}

After identifying $\mathbb C$ with ${\mathbb R}^2$ or ${\mathbb P}^2_+$,
the first observation is that numerical range ${\mathcal W}(A)$
is dual to LMI set ${\mathcal F}(A)$, and hence it is 
an affine projection of the semidefinite cone.

\begin{lemma}\label{proj}
${\mathcal W}(A) = {\mathcal F}(A)^* = \{\left(\langle A_0,W \rangle,\:
\langle A_1,W \rangle,\:\langle A_2,W \rangle\right) \in {\mathbb P}^2_+
\: :\: W \in {\mathbb C}^{n\times n},\: W \succeq 0\}$.
\end{lemma}

{\bf Proof}:
The dual to ${\mathcal F}(A)$ is
\[
\begin{array}{rcl}
{\mathcal F}(A)^* & = & \{x \: :\: \langle x,y \rangle =
\langle F(y),W \rangle = \sum_k \langle A_k,W \rangle  y_k 
\geq 0 , \: W \succeq 0\} \\
& = & \{x \: :\: x_k = \langle A_k,W \rangle, \: 
W \succeq 0\},
\end{array}
\]
an affine projection of the semidefinite cone.
On the other hand, since $w^*Aw = w^*A_1w + i\: w^*A_2w$,
the numerical range can be expressed as
\[
\begin{array}{rcl}
{\mathcal W}(A) & = & \{x = (w^*A_0w,\:w^*A_1w,\:w^*A_2w)\} \\
& = & \{x \: :\: x_k = \langle A_k,W \rangle, \:
W \succeq 0,\:\mathrm{rank}\:W = 1\},
\end{array}
\]
the same affine projection as above, acting now on a subset
of the semidefinite cone, namely the non-convex variety of
rank-one positive semidefinite matrices $W = ww^*$. Since $w^*A_0w = 1$,
set ${\mathcal W}(A)$ is compact, and $\mathrm{conv}\:{\mathcal W}(A) =
{\mathcal F}(A)^*$. The equality ${\mathcal W}(A) = {\mathcal F}(A)^*$
follows from the Toeplitz-Hausdorff theorem establishing convexity
of ${\mathcal W}(A)$, see \cite[Section 1.3]{horn} or
\cite[Theorem 1.1-2]{gustafson}.
$\Box$

Lemma \ref{proj} indicates that the numerical range
has the geometry of planar projections of
the semidefinite cone. In the terminology 
of \cite[Chapter 4]{bental},
the numerical range is semidefinite representable.

\section{Convex hulls of algebraic curves}

In this section, we notice that the boundaries of numerical range
${\mathcal W}(A)$ and its dual LMI set ${\mathcal F}(A)$ are subsets of
algebraic curves $\mathcal P$ and $\mathcal Q$ defined respectively
in (\ref{curvep}) and (\ref{curveq}), and explicitly given as locii
of determinants of Hermitian pencils.

\subsection{Dual curve}

\begin{lemma}\label{fap}
${\mathcal F}(A)$ is the connected component
delimited by $\mathcal P$ around the origin.
\end{lemma}

{\bf Proof}: A ray starting from the origin leaves LMI set ${\mathcal F}(A)$
when the determinant $p(y) = \det\:\sum_k y_k A_k$ vanishes. Therefore
the boundary of ${\mathcal F}(A)$ is the subset of algebraic curve $\mathcal P$
belonging to the convex connected component containing the origin.$\Box$

Note that $\mathcal P$, by definition, is the locus, or vanishing set of
a determinant of a Hermitian pencil. Moreover, the pencil is definite
at the origin so the corresponding polynomial $p(y)$ satisfies
a real zero (hyperbolicity) condition. Connected components delimited by such
determinantal locii are studied in \cite{helton}, where it is shown
that they correspond to feasible sets of two-dimensional LMIs. A remarkable
result of \cite{helton} is that every planar LMI set can be expressed this way.
These LMI sets form a strict subset of planar convex basic semi-algebraic sets, called
rigidly convex sets (see \cite{helton} for examples of convex basic
semi-algebraic sets which are not rigidly convex).
Rigidly convex sets are affine sections of the semidefinite cone.

\subsection{Primal curve}

\begin{lemma}\label{waq}
${\mathcal W}(A) = \mathrm{conv}\:{\mathcal Q}$.
\end{lemma}

{\bf Proof}: From the proof of Lemma \ref{proj},
a supporting line $\{x \: :\: \sum_k x_k y_k = 0\}$
to ${\mathcal W}(A)$ has coefficients $y$ satisfying $p(y)=0$. 
The boundary of ${\mathcal W}(A)$ is therefore generated
as an envelope of the supporting lines.
See \cite{murnaghan}, \cite[Theorem 10]{kippenhahn} and also
\cite[Theorem 1.3]{fiedler}.$\Box$

${\mathcal Q}$ is called the boundary generating curve of matrix $A$
in \cite{kippenhahn}. An interesting feature is that, similarly to
$\mathcal P$, curve $\mathcal Q$ can
be expressed as the locus of a determinant of a Hermitian pencil.
In the case ${\mathcal Q}$ is irreducible (i.e. polynomial $q(x)$
cannot be factored) and ${\mathcal P}$ is not singular (i.e. there
is no point in the complex projective plane
such that the gradient of $p(x)$ vanishes) then $q(x)$ can be
written (up to a multiplicative constant) as the determinant
of a symmetric pencil, see \cite[Theorem 2.4]{fiedler}.
Discrete differentials and B\'ezoutians can also be used to
construct symmetric affine determinantal representations,
see \cite[Section 4.2]{henrion}. Note however that the constructed
pencils are not sign definite.
Hence the convex hull ${\mathcal W}(A)$ is not
a rigidly convex LMI set, it cannot be an affine section
of the semidefinite cone. However, as noticed in Lemma \ref{proj},
it is an affine projection of the semidefinite cone.

\section{Examples}

\subsection{Rational cubic and quartic}

\begin{figure}[h!]
\begin{minipage}[t]{8cm}
\begin{center}
\includegraphics[width=8cm]{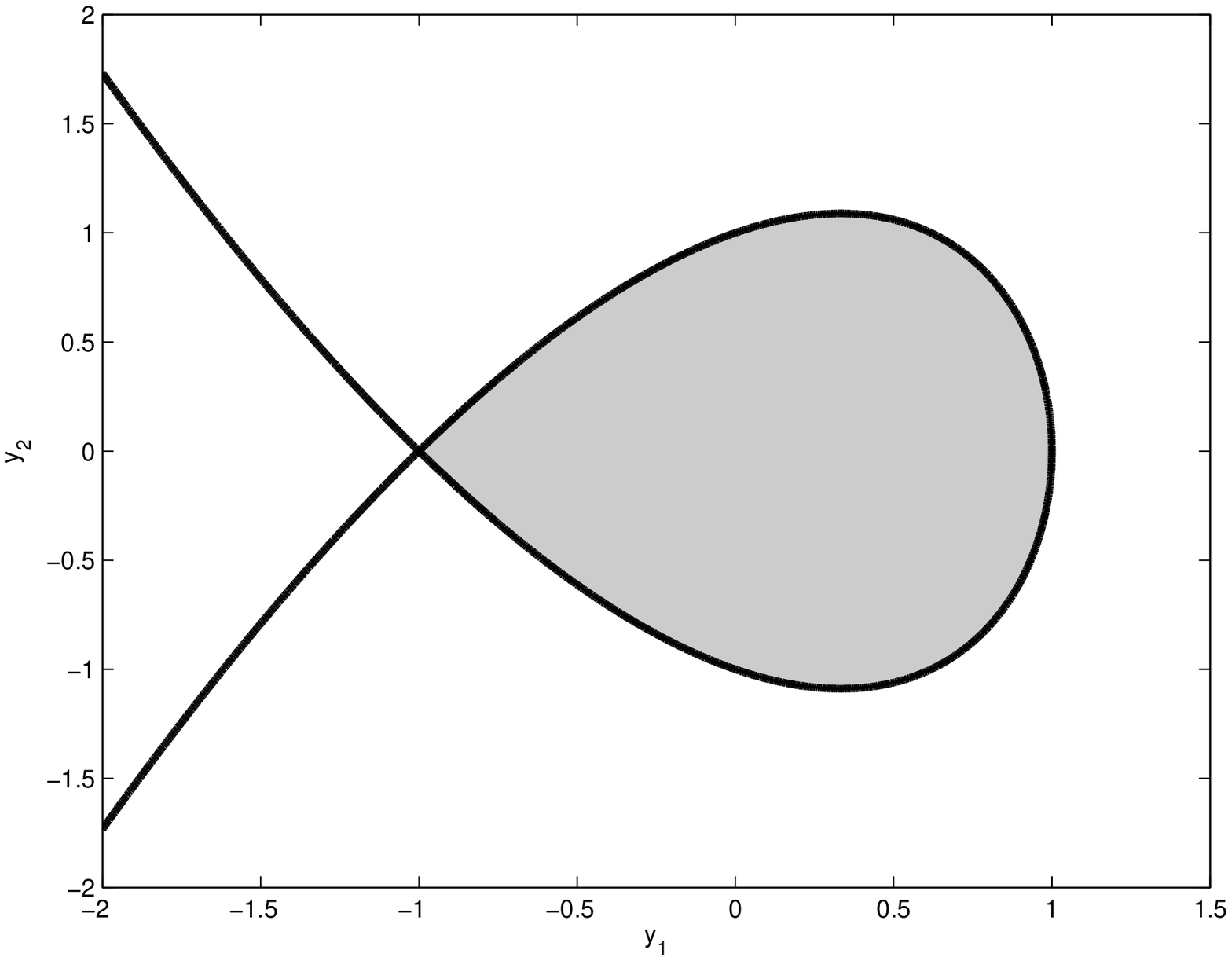}
\end{center}
\end{minipage}
\begin{minipage}[t]{8cm}
\begin{center}
\includegraphics[width=8cm]{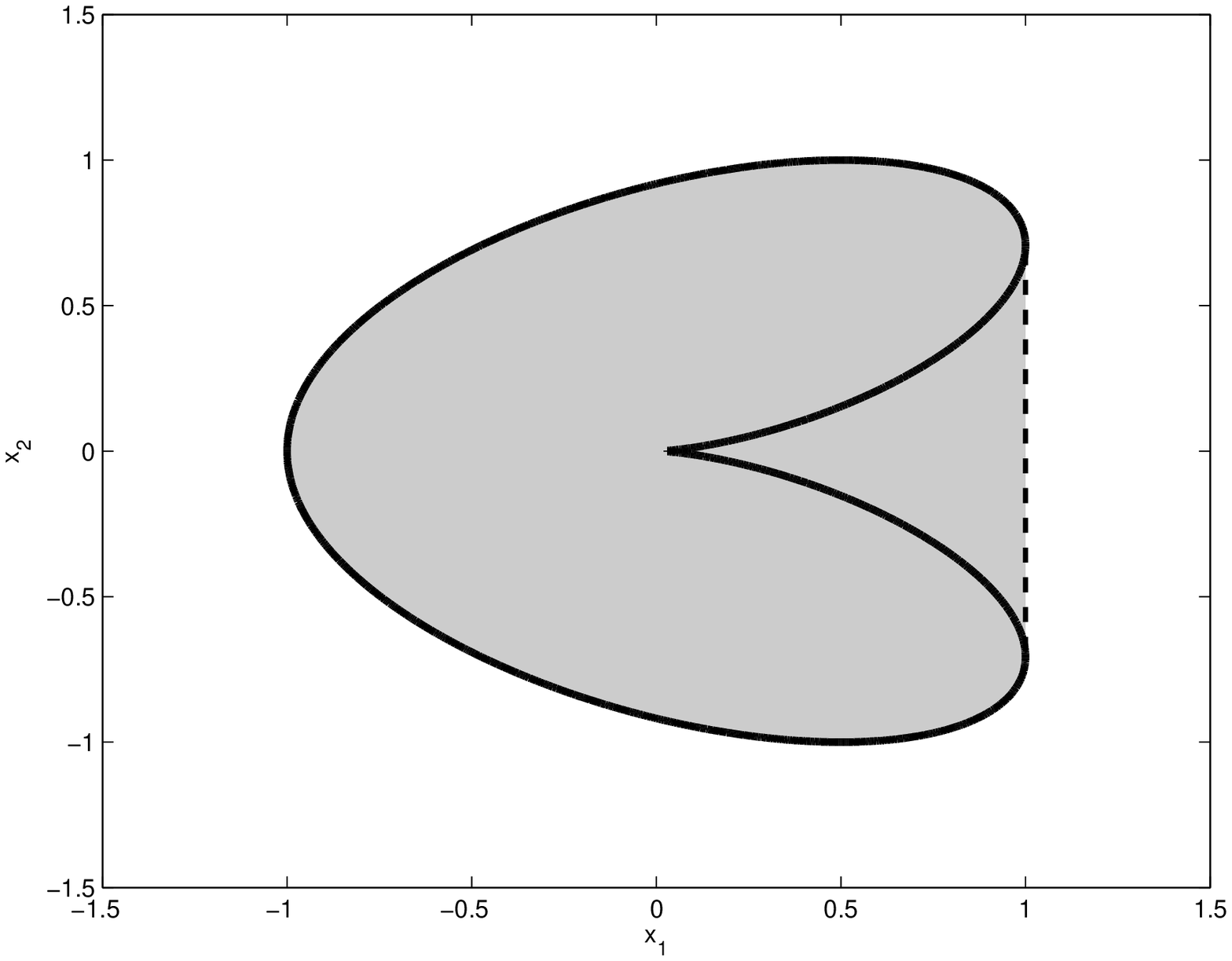}
\end{center}
\end{minipage}
\caption{Left: LMI set ${\mathcal F}(A)$ (gray area) delimited by
cubic $\mathcal P$ (black). Right: numerical range ${\mathcal W}(A)$
(gray area, dashed line) convex hull
of quartic $\mathcal Q$ (black solid line).\label{ex1}}
\end{figure}

Let
\[
A = \left[\begin{array}{ccc}
0 & 0 & 1 \\ 0 & 1 & i \\ 1 & i & 0 
\end{array}\right].
\]
Then
\[
F(y) = \left[\begin{array}{ccc}
y_0 & 0 & y_1 \\
0 & y_0 + y_1 & y_2 \\
y_1 & y_2 & y_0
\end{array}\right]
\]
and
\[
p(y) = (y_0-y_1)(y_0+y_1)^2-y_0y^2_2
\]
defines a genus-zero cubic curve $\mathcal P$
whose connected component containing the origin
is the LMI set ${\mathcal F}(A)$, see Figure \ref{ex1}.
With an elimination technique (resultants or Gr\"obner basis with
lexicographical ordering), we obtain
\[
q(x) = 4x^4_1+32x^4_2+13x^2_1x^2_2-18x_0x_1x^2_2+4x_0x^3_1-27x^2_0x^2_2
\]
defining the dual curve $\mathcal Q$, a genus-zero quartic
with a cusp, whose convex hull is the numerical range 
${\mathcal W}(A)$, see Figure \ref{ex1}.

\subsection{Couple of two nested ovals}

\begin{figure}[h!]
\begin{minipage}[t]{8cm}
\begin{center}
\includegraphics[width=8cm]{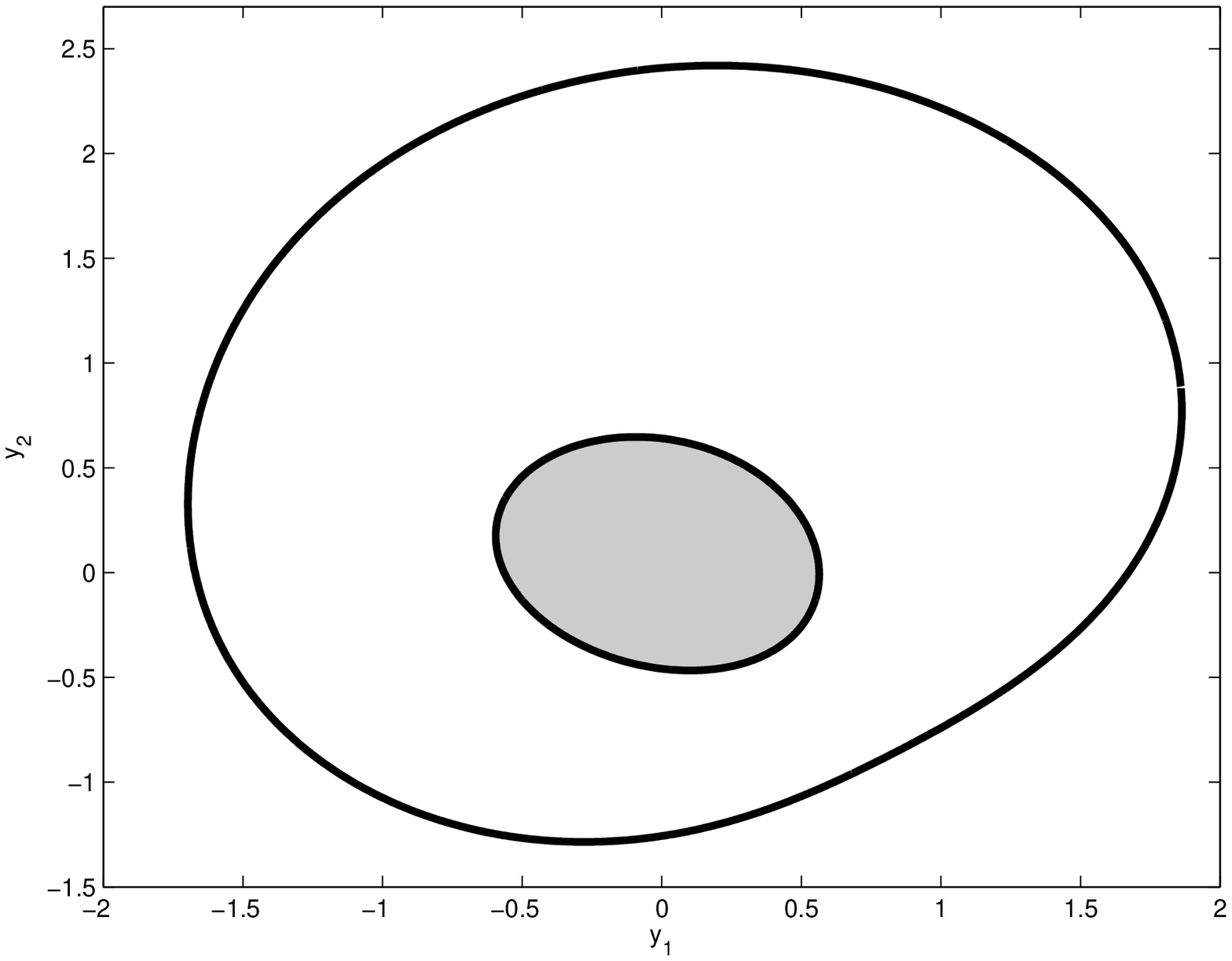}
\end{center}
\end{minipage}
\begin{minipage}[t]{8cm}
\begin{center}
\includegraphics[width=8cm]{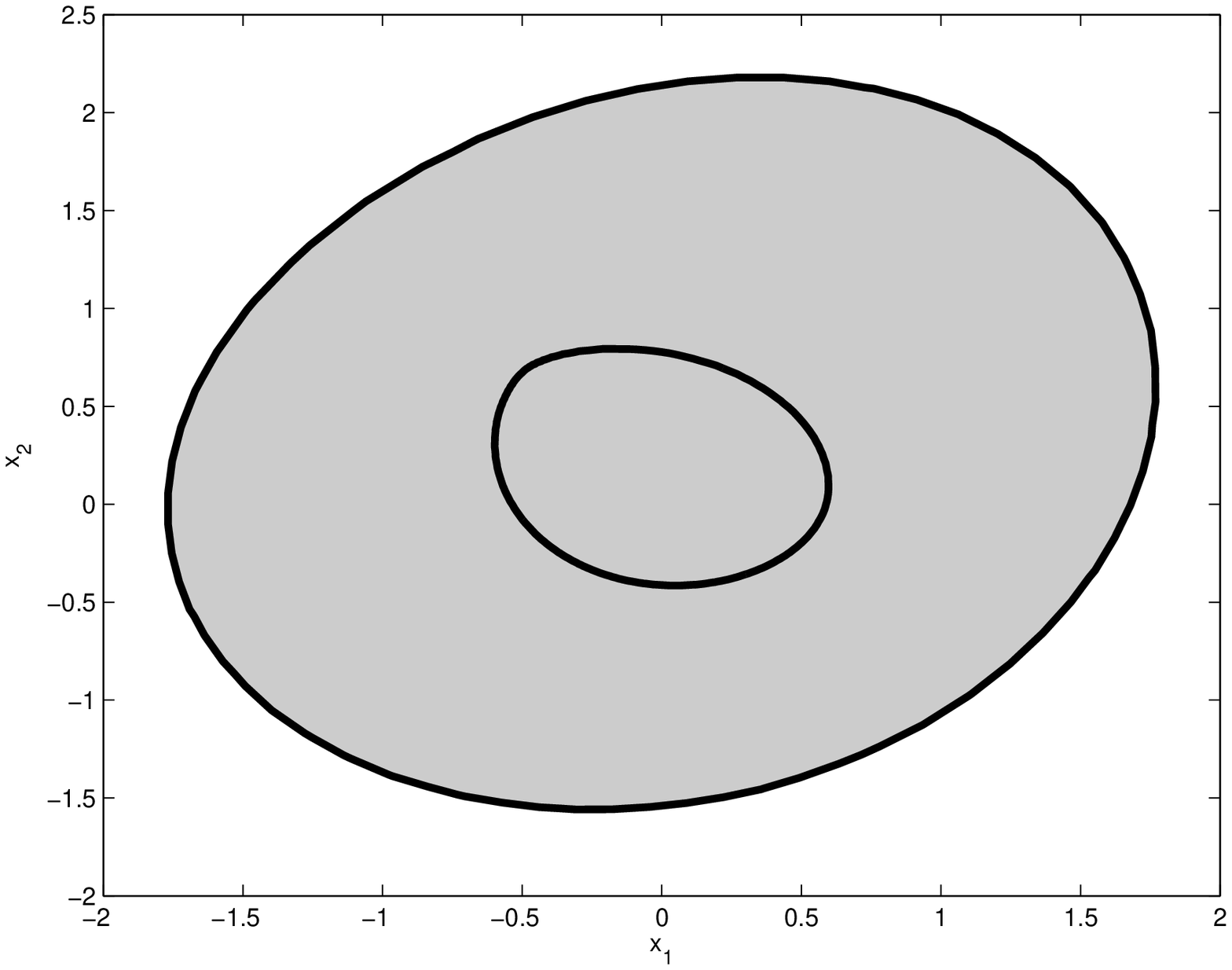}
\end{center}
\end{minipage}
\caption{Left: LMI set ${\mathcal F}(A)$ (gray area) delimited by the
inner oval of quartic $\mathcal P$ (black line). Right: numerical range ${\mathcal W}(A)$
(gray area) delimited by the outer oval of octic $\mathcal Q$ (black line).\label{ex2}}
\end{figure}

For
\[
A = \left[\begin{array}{cccc}
0 & 2 & 1+2i & 0 \\ 0 & 0 & 1 & 0 \\ 0 & i & i & 0 \\ 0 & -1+i & i & 0 
\end{array}\right]
\]
the quartic $\mathcal P$ and its dual octic $\mathcal Q$ both feature
two nested ovals, see Figure \ref{ex2}. The inner oval delimited by
$\mathcal P$ is rigidly convex, whereas the outer oval delimited
by $\mathcal Q$ is convex, but not rigidly convex.

\subsection{Cross and star}

\begin{figure}[h!]
\begin{minipage}[t]{8cm}
\begin{center}
\includegraphics[width=8cm]{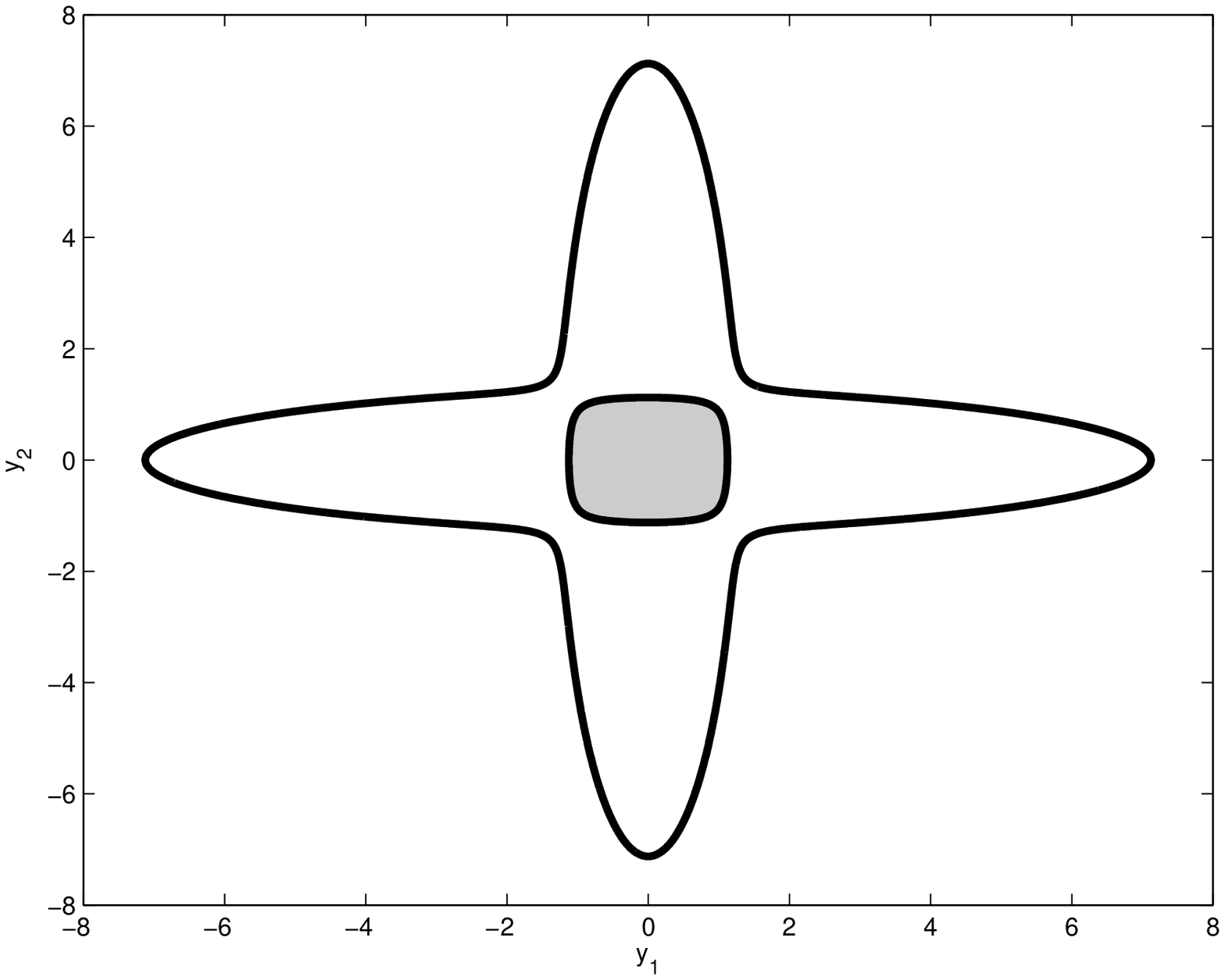}
\end{center}
\end{minipage}
\begin{minipage}[t]{8cm}
\begin{center}
\includegraphics[width=8cm]{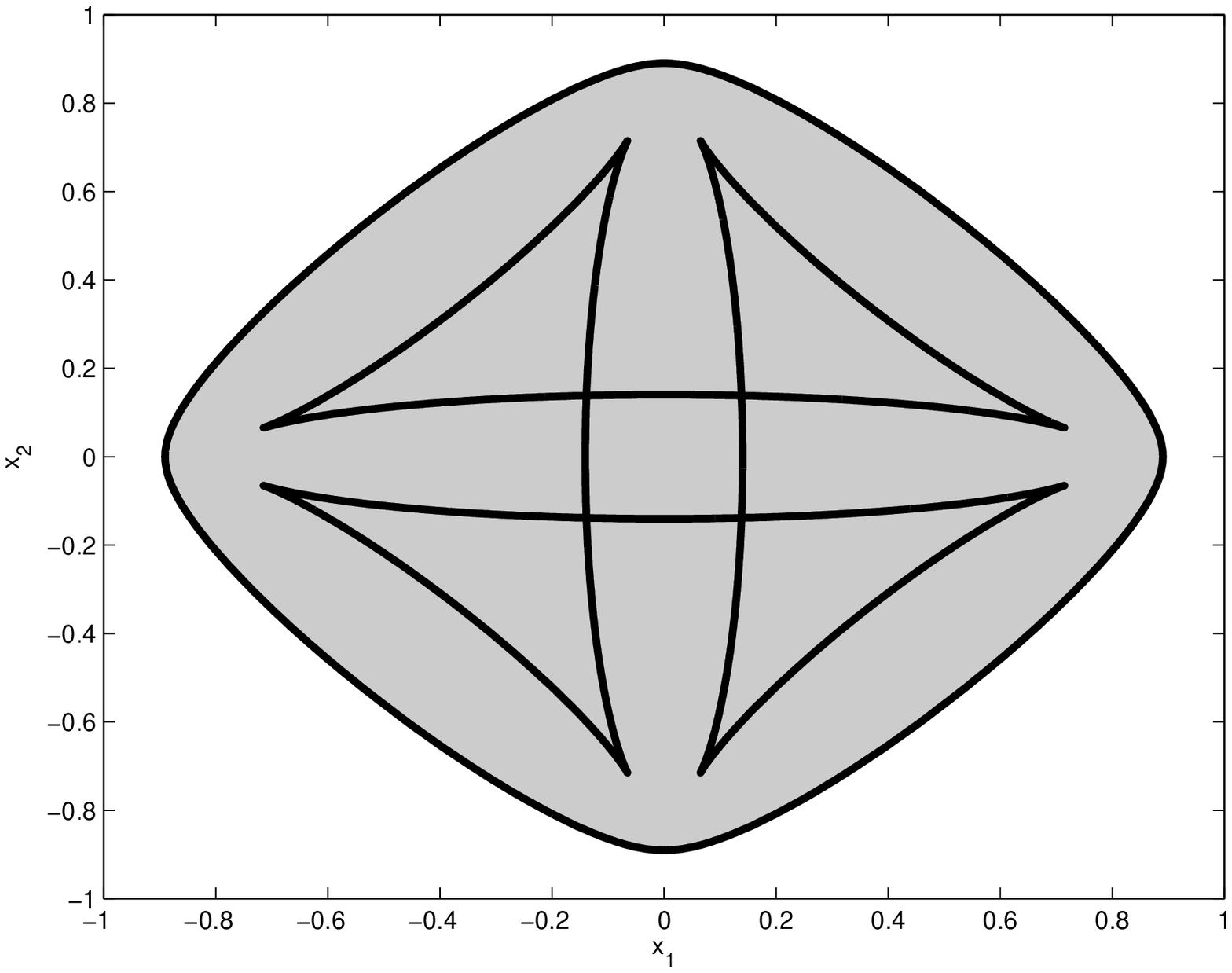}
\end{center}
\end{minipage}
\caption{Left: LMI set ${\mathcal F}(A)$ (gray area) delimited by the
inner oval of quartic $\mathcal P$ (black line). Right: numerical range ${\mathcal W}(A)$
(gray area) delimited by the outer oval of twelfth-degree $\mathcal Q$ (black line).\label{gustafson1}}
\end{figure}

A computer-generated representation of the numerical range as
an enveloppe curve can be found in \cite[Figure 1, p. 139]{gustafson}
for
\[
A = \left[\begin{array}{cccc}
0 & 1 & 0 & 0 \\ 0 & 0 & 1 & 0 \\ 0 & 0 & 0 & 1 \\ \frac{1}{2} & 0 & 0 & 0 
\end{array}\right].
\]
We obtain the quartic
\[
p(y) = \frac{1}{64}(64 y_0^4 - 52 y_0^2 y_1^2 - 52 y_0^2 y_2^2 
+ y_1^4 + 34 y_1^2 y_2^2 + y_2^4)
\]
and the dual twelfth-degree polynomial
\[
\begin{array}{rcl}
q(x) & = & 5184x_0^{12}-299520x_0^{10}x_1^2-299520x_0^{10}x_2^2+1954576x_0^8x_1^4\\
&& +16356256x_0^8x_1^2x_2^2+1954576x_0^8x_2^4-5375968x_0^6x_1^6-79163552x_0^6x_1^4x_2^2\\
&&-79163552x_0^6x_1^2x_2^4-5375968x_0^6x_2^6+7512049x_0^4x_1^8+152829956x_0^4x_1^6x_2^2\\
&&-2714586x_0^4x_1^4x_2^4+152829956x_0^4x_1^2x_2^6+7512049x_0^4x_2^8-5290740x_0^2x_1^{10}\\
&&-136066372x_0^2x_1^8x_2^2+232523512x_0^2x_1^6x_2^4+232523512x_0^2x_1^4x_2^6-136066372x_0^2x_1^2x_2^8\\
&&-5290740x_0^2x_2^{10}+1498176x_1^{12}+46903680x_1^{10}x_2^2-129955904x_1^8x_2^4\\
&&+186148096x_1^6x_2^6-129955904x_1^4x_2^8+46903680x_1^2x_2^{10}+1498176x_2^{12}
\end{array}
\]
whose corresponding curves and convex hulls are represented in Figure \ref{gustafson1}.

\subsection{Decomposition into irreducible factors}

\begin{figure}[h!]
\begin{minipage}[t]{8cm}
\begin{center}
\includegraphics[width=8cm]{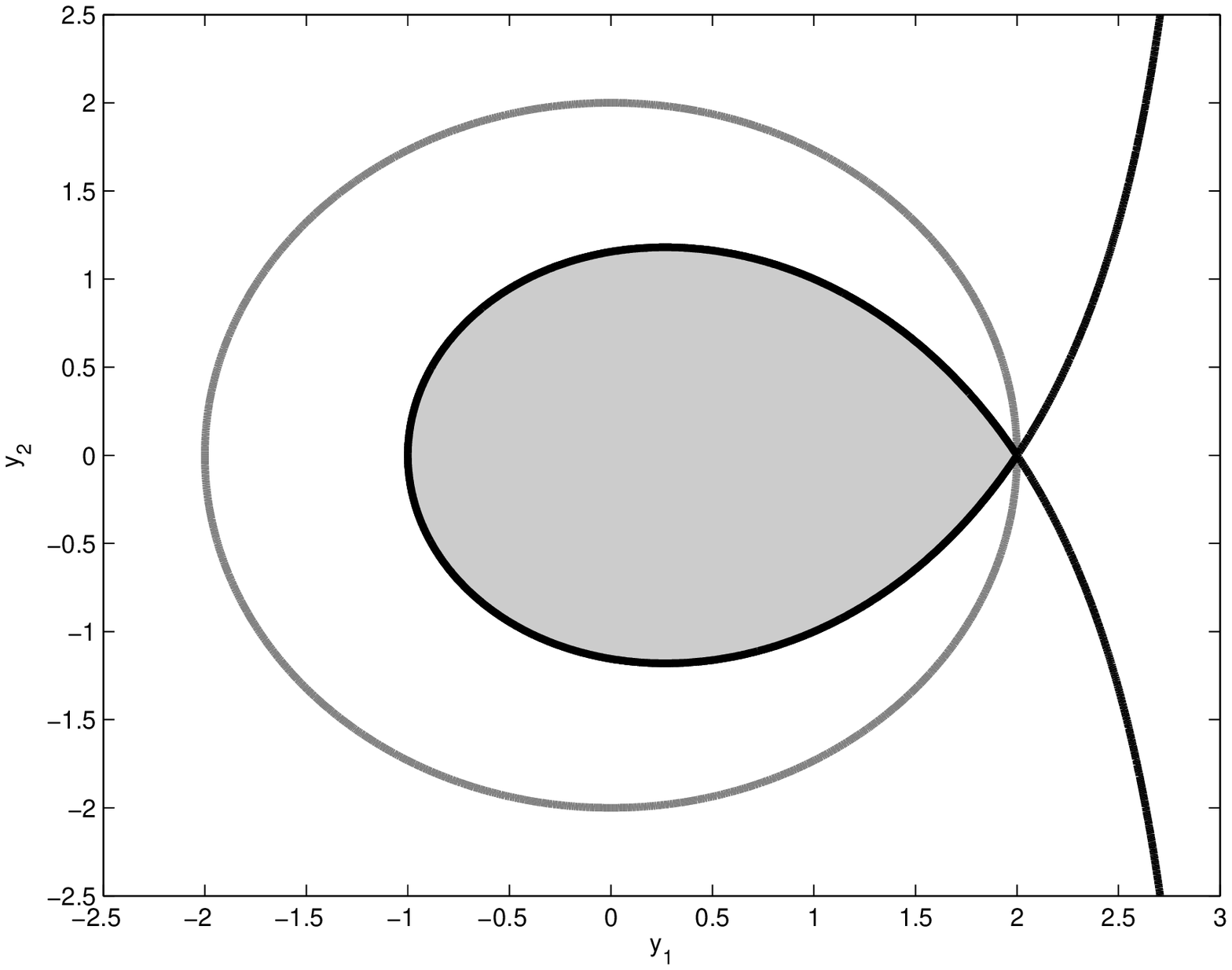}
\end{center}
\end{minipage}
\begin{minipage}[t]{8cm}
\begin{center}
\includegraphics[width=8cm]{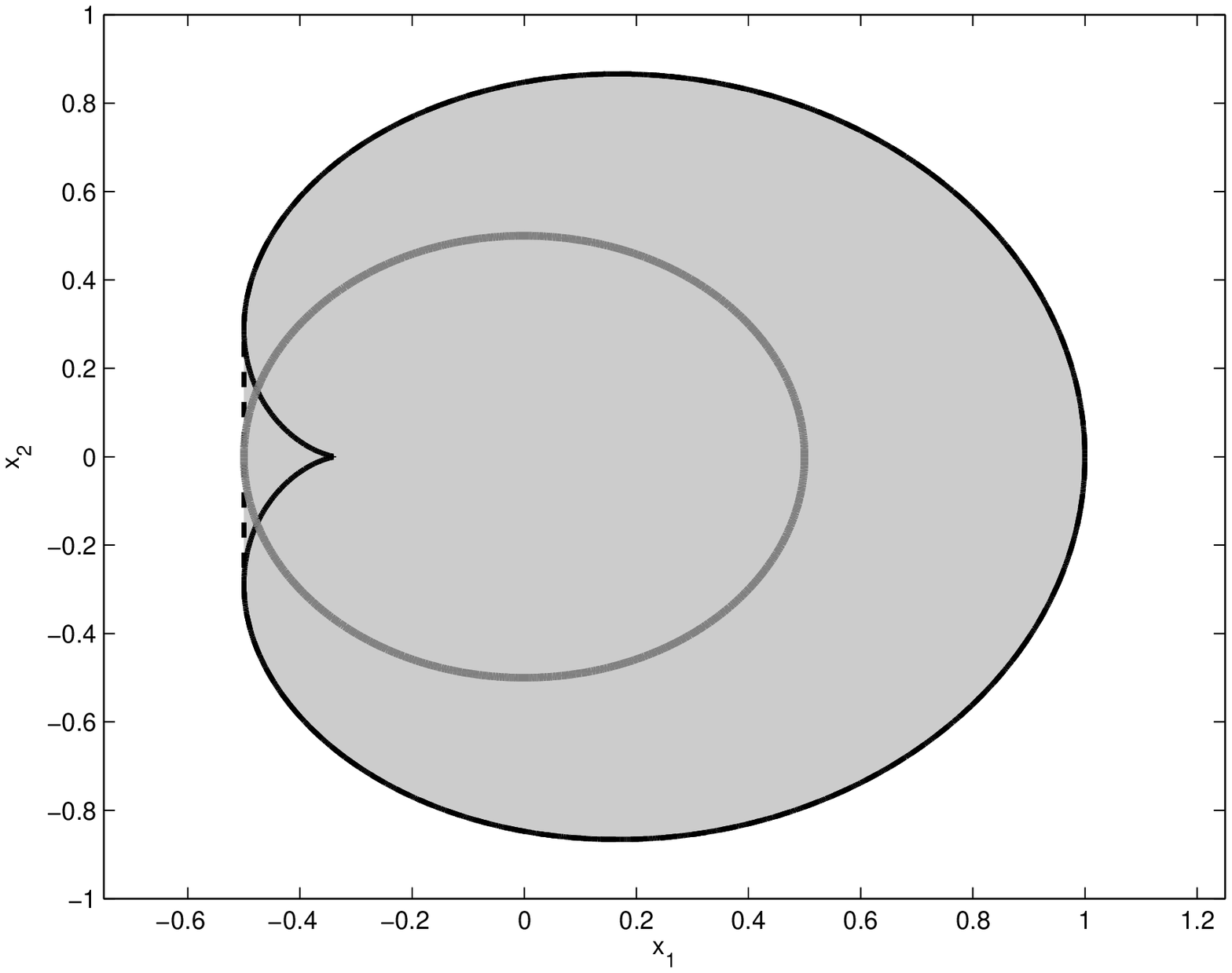}
\end{center}
\end{minipage}
\caption{Left: LMI set ${\mathcal F}(A)$ (gray area) intersection of cubic (black solid line)
and conic (gray line) LMI sets. Right: numerical range ${\mathcal W}(A)$
(gray area, black dashed line) convex hull of the union of a quartic curve (black solid line)
and conic curve (gray line).\label{gustafson6}}
\end{figure}

Consider the example of \cite[Figure 6, p. 144]{gustafson} with
\[
A = 
\left[\begin{array}{ccccccccc}
0 & 0 & 0 & 0 & 1 & 0 & 0 & 0 & 1 \\
0 & 0 & 0 & 0 & 0 & 1 & 0 & 0 & 0 \\
0 & 0 & 0 & 0 & 0 & 0 & 1 & 0 & 0 \\
0 & 0 & 0 & 0 & 0 & 0 & 0 & 1 & 0 \\
0 & 0 & 0 & 0 & 0 & 0 & 0 & 0 & 1 \\
0 & 0 & 0 & 0 & 0 & 0 & 0 & 0 & 0 \\
0 & 0 & 0 & 0 & 0 & 0 & 0 & 0 & 0 \\
0 & 0 & 0 & 0 & 0 & 0 & 0 & 0 & 0 \\
0 & 0 & 0 & 0 & 0 & 0 & 0 & 0 & 0 \\
\end{array}\right]
\]
The determinant of the trivariate pencil
factors as follows
\[
p(y) = \frac{1}{256}(4y_0^3-3y_0y_1^2-3y_0y_2^2+y_1^3+y_1y_2^2)(4y_0^2-y_1^2-y_2^2)^3
\]
which means that the LMI set ${\mathcal F}(A)$ is the intersection
of a cubic and conic LMI. 

The dual curve $\mathcal Q$ is the union of the quartic
\[
{\mathcal Q}_{\;1}=\{x \: :\: x_0^4-8x_0^3x_1-18x_0^2x_1^2 -18x_0^2x_2^2+27x_1^4+54x_1^2x_2^2+27x_2^4=0\},
\]
a cardioid dual to the cubic factor of $p(y)$, and the conic
\[
{\mathcal Q}_{\;2}=\{x \: :\: x_0^2-4x_1^2-4x_2^2=0\},
\]
a circle dual to the quadratic factor of $p(y)$.
The numerical range ${\mathcal W}(A)$ is the convex hull of the union
of $\mathrm{conv}\:{\mathcal Q}_{\;1}$ and $\mathrm{conv}\:{\mathcal Q}_{\;2}$,
which is here the same as $\mathrm{conv}\:{\mathcal Q}_{\;1}$, see Figure \ref{gustafson6}.

\subsection{Polytope}\label{polytope}

\begin{figure}[h!]
\begin{minipage}[t]{8cm}
\begin{center}
\includegraphics[width=8cm]{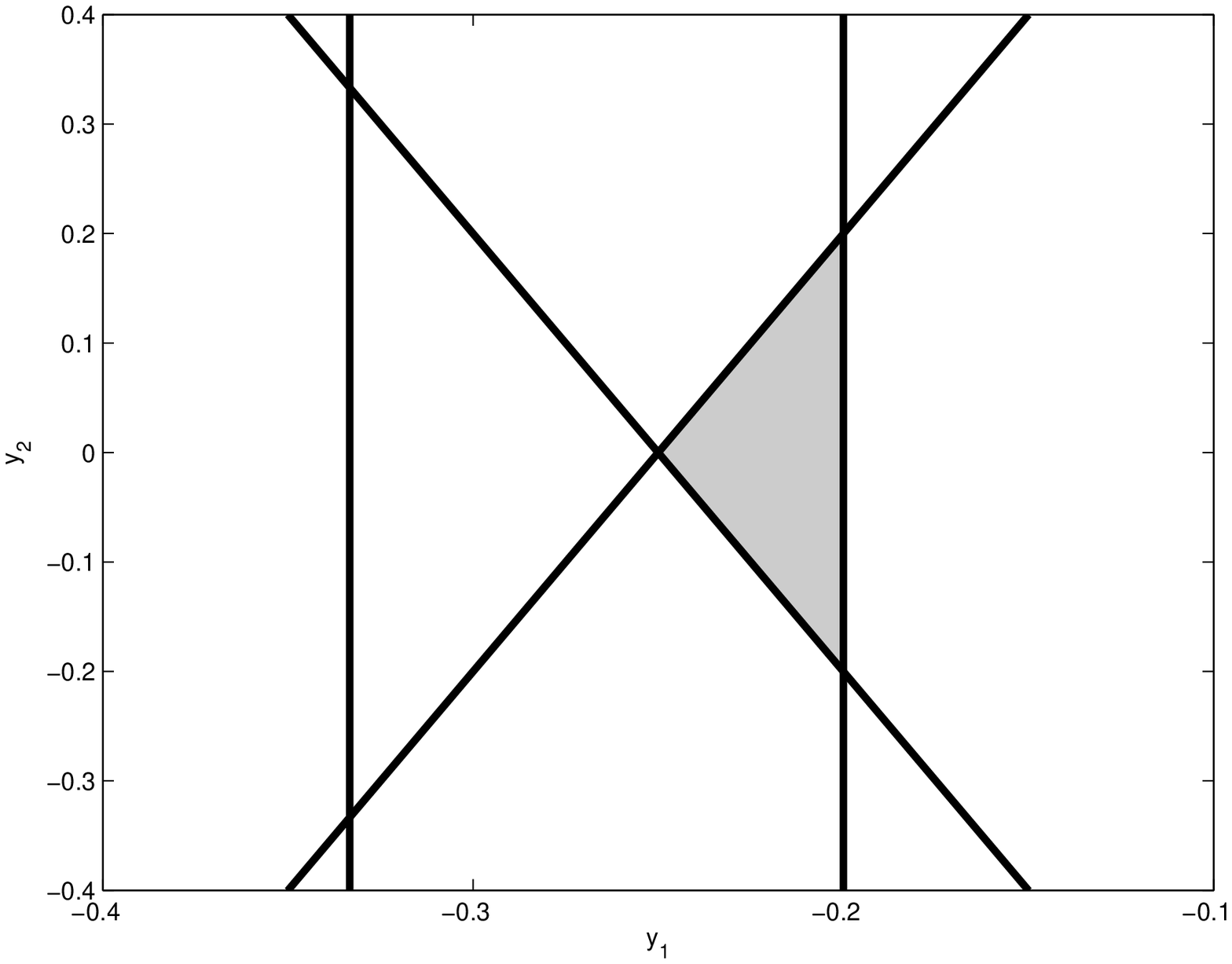}
\end{center}
\end{minipage}
\begin{minipage}[t]{8cm}
\begin{center}
\includegraphics[width=8cm]{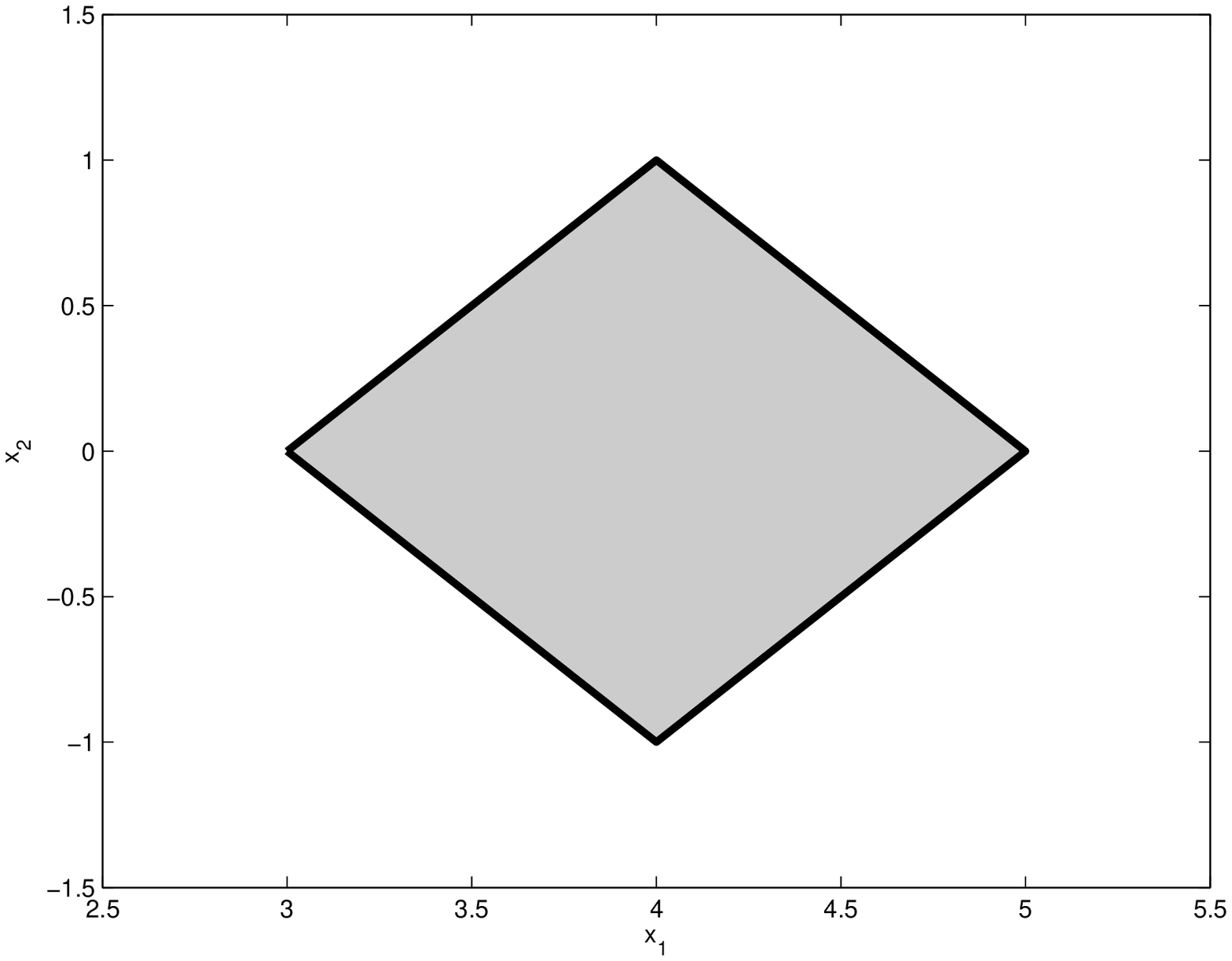}
\end{center}
\end{minipage}
\caption{Left: LMI set ${\mathcal F}(A)$ (gray area) intersection of four half-planes.
Right: numerical range ${\mathcal W}(A)$
(gray area) convex hull of four vertices.\label{gustafson9}}
\end{figure}

Consider the example of \cite[Figure 9, p. 147]{gustafson} with
\[
A = 
\left[\begin{array}{rrrr}
4 & 0 & 0 & -1 \\
-1 & 4 & 0 & 0 \\
0 & -1 & 4 & 0 \\
0 & 0 & -1 & 4
\end{array}\right].
\]
The dual determinant factors into linear terms
\[
p(y) = (y_0+5y_1)(y_0+3y_1)(y_0+4y_1+y_2)(y_0+4y_1-y_2)
\]
and this generates a polytopic LMI set ${\mathcal F}(A) = \{y \: :\:
y_0+5y_1 \geq 0, \:y_0+3y_1 \geq 0, y_0+4y_1+y_2 \geq 0,
y_0+4y_1-y_2 \geq 0\}$, a triangle with vertices
$(1,-\frac{1}{4},0)$, $(1,-\frac{1}{5},\frac{1}{5})$
and $(1,-\frac{1}{5},-\frac{1}{5})$. The dual to curve $\mathcal P$
is the union of the four points $(1,5,0)$, $(1,3,0)$,
$(1,4,1)$ and $(1,4,-1)$ and hence the numerical range
${\mathcal W}(A)$ is the polytopic convex hull of these four vertices,
see Figure \ref{gustafson9}.

\section{A problem in statistics}

We have seen with Example \ref{polytope} that the numerical range
can be polytopic, and this is the case in particular when $A$
is a normal matrix (i.e. satisfying $A^*A=AA^*$), see e.g.
\cite[Theorem 3]{kippenhahn} or
\cite[Theorem 1.4-4]{gustafson}.

In this section, we study a problem
that boils down to studying rectangular numerical ranges, i.e.
polytopes with edges parallel to the main axes.
Craig's theorem is a result
from statistics on the stochastic independence of two quadratic
forms in variates following a joint normal distribution,
see \cite{driscoll} for 
an historical account. In its simplest form (called the
central case) the result can be stated as follows (in the
sequel we work in the affine plane $y_0=1$):

\begin{theorem}\label{det}
Let $A_1$ and $A_2$ be Hermitian matrices of size $n$. Then
$\det(I_n+y_1A_1+y_2A_2) = \det(I_n+y_1A_1)\det(I_n+y_2A_2)$
if and only if $A_1 A_2=0$.
\end{theorem}

{\bf Proof}:
If $A_1 A_2 = 0$ then obviously $\det(I_n+y_1A_1)\det(I_n+y_2A_2)
=\det((I_n+y_1A_1)(I_n+y_2A_2))=\det(I_n+y_1A_1+y_2A_2+y_1y_2A_1A_2)=
\det(I_n+y_1A_1+y_2A_2)$. Let us prove the converse statement.

Let ${a_1}_k$ and ${a_2}_k$ respectively denote the eigenvalues of 
$A_1$ and $A_2$, for $k=1,\ldots,n$. Then $p(y)=\det(I_n+y_1A_1+y_2A_2)
=\det(I_n+y_1A_1)\det(I_n+y_2A_2)=\prod_k(1+y_1{a_1}_k)\prod_k(1+y_2{a_2}_k)$
factors into linear
terms, and we can write $p(y)=\prod_k(1+y_1{a_1}_k+y_2{a_2}_k)$
with ${a_1}_k{a_2}_k=0$ for all $k=1,\ldots,n$.
Geometrically, this means that the corresponding numerical range
${\mathcal W}(A)$ for $A=A_1+iA_2$
is a rectangle with vertices $(\min_k a_k,\min_k b_k)$,
$(\min_k a_k,\max_k b_k)$, $(\max_k a_k,\min_k b_k)$ and
$(\max_k a_k,\max_k b_k)$. 

Following the terminology of \cite{motzkin}, $A_1$ and $A_2$
satisfy property L since $y_1A_1+y_2A_2$ has
eigenvalues $y_1 {a_1}_k + y_2 {a_2}_k$ for $k=1,\ldots,n$.
From \cite[Theorem 2]{motzkin} it
follows that $A_1 A_2 = A_2 A_1$, and hence that the two
matrices are simultaneously diagonalisable: there exists
a unitary matrix $U$ such that $U^*A_1U = \mathrm{diag}_k {a_1}_k$
and $U^*A_2U = \mathrm{diag}_k {a_2}_k$. Since ${a_1}_k {a_2}_k=0$
for all $k$, we have $\sum_k {a_1}_k{a_2}_k = U^*A_1UU^*A_2U = U^*A_1A_2U = 0$
and hence $A_1A_2=0$.
$\Box$

\section{Conclusion}

The geometry of the numerical range, studied to a large extent
by Kippenhahn in \cite{kippenhahn} -- see \cite{zachlin} for
an English translation with comments and corrections -- is
revisited here from the perspective of semidefinite
programming duality. In contrast with previous studies of
the geometry of the numerical range, based on differential
topology \cite{jonckheere}, it is namely noticed that the numerical
range is a semidefinite representable set, an affine projection
of the semidefinite cone, whereas its geometric dual is
an LMI set, an affine section of the semidefinite cone. The boundaries
of both primal and dual sets are components of algebraic
plane curves explicitly formulated as locii
of determinants of Hermitian pencils. The geometry of
the numerical range is therefore the geometry of (planar sections
and projections of) the semidefinite cone, and hence every study
of this cone is also relevant to the study of the numerical range.

The notion of numerical range can be generalized in various directions,
for example in spaces of dimension greater than two, where
it is non-convex in general \cite{fan}. Its convex hull is
still representable as a projection of the semidefinite cone,
and this was used extensively in the scope of robust
control to derive computationally tractable but potentially
conservative LMI stability conditions for uncertain
linear systems, see e.g. \cite{packard}.
The numerical range of three matrices is mentioned in
\cite[Section 1.8]{horn}. In this context,
it would be interesting to derive conditions on three
matrices $A_1$, $A_2$, $A_3$ ensuring
that $\det(I_n+y_1A_1+y_2A_2+y_3A_3) = \det(I_n+y_1A_1)
\det(I_n+y_2A_2)\det(I_n+y_3A_3)$.
Another extension of
the numerical range to matrix polynomials (including matrix
Pencils) was carried out in \cite{chien}, also using
algebraic geometric considerations, and it could be
interesting to study semidefinite representations of
convex hulls of these numerical ranges.

The inverse problem of finding a matrix given its numerical
range (as the convex hull of a given algebraic curve)
seems to be difficult. In a sense, it is dual to the
problem of finding a symmetric (or Hermitian) definite linear
determinantal representation of a trivariate form:
given $p(y)$ satisfying a real zero (hyperbolicity)
condition, find Hermitian matrices $A_k$ such that
$p(y)=\det(\sum_k y_k A_k)$, with $A_0$ positive definite.
Explicit formulas are
described in \cite{helton} based on transcendental
theta functions and Riemann surface theory, and
the case of curves $\{y \: :\: p(y)=0\}$ of genus zero is settled in \cite{henrion}
using B\'ezoutians. A more direct and computationally viable
approach in the positive genus case is still missing,
and one may wonder whether the geometry of the dual
object, namely the numerical range $\mathrm{conv}\{x \: :\: q(x)=0\}$,
could help in this context.

\section*{Acknowledgments}

The author is grateful to Leiba Rodman for his suggestion of
studying rigid convexity of the numerical range. This work
also benefited from technical advice by Jean-Baptiste Hiriart-Urruty
who recalled Theorem \ref{det} in the September 2007 issue of the MODE
newsletter of SMAI (French society for applied and industrial mathematics)
and provided reference \cite{driscoll}. Finally,
references \cite{chien,jonckheere} were brought to the author's
attention by an anonymous referee.

\end{document}